\documentclass{elsart}
\usepackage{amssymb}
\usepackage{latexsym,amsmath,pstricks,pst-node, amsfonts, epsfig}
\newcommand{\Aut}{\mathop{\mathrm{Aut}}}
\newcommand{\Cay}{\mathop{\mathrm{Cay}}}

\newtheorem{theorem}{Theorem}[section]
\newtheorem{lemma}[theorem]{Lemma}
\newtheorem{sublemma}{Claim}[theorem]
\newtheorem{definition}[theorem]{Definition}

\newtheorem{conjecture}[theorem]{Conjecture}
\numberwithin{figure}{section}
\begin{document}
\begin{frontmatter}
\title{Balanced Cayley graphs and balanced planar graphs}
\author{Joy Morris\thanksref{joy}}, \author{Pablo Spiga\thanksref{mouse}} and \author{Kerri Webb}
\thanks[joy]{This research was supported by NSERC of Canada}
\thanks[mouse]{This research was supported by a fellowship
from the Pacific Institute for the Mathematical Sciences\\
\emph{Email address:} \texttt{spiga@cs.uleth.ca} (Pablo Spiga)}
\address{University of Lethbridge\\Department of Mathematics and
Computer Science\\ $4401$ University Drive, Lethbridge, Alberta,
Canada}
\begin{abstract} A balanced graph
is a bipartite graph with no induced circuit of length $2\pmod 4$.
These graphs arise in linear programming. We focus on
graph-algebraic properties of balanced graphs to prove a complete
classification of balanced Cayley graphs on abelian groups.
Moreover, in Section~\ref{sec:planar} of this paper, we prove that
there is no cubic balanced planar graph. Finally, some remarkable
conjectures for balanced regular graphs are also presented.
\end{abstract}
\end{frontmatter}

\section{Introduction}
A $\{0,1\}$-matrix is \emph{balanced} if the sum of the entries of
every submatrix that is minimal with respect to the property of
containing $2$ nonzero entries per row and per column, is congruent to
$0\pmod 4$. Balanced matrices were introduced by
Berge~\cite{ber:bal} in the context of hypergraphs, and they arise
naturally in linear programming~\cite{Ta}. 

There has been considerable study of
balanced matrices; the reader might
check~\cite{con:fag} or~\cite{cor:com} for a survey on the main
results and horizons on balanced matrices.

Every $\{0,1\}$-matrix is also the bipartite adjacency matrix of a bipartite
graph. Specifically, if $X$ is a bipartite graph with vertex bipartition
$(U,V)$, then the bipartite adjacency matrix for $X$ is the
$\{0,1\}$-matrix $A$ with $A_{u,v}=1$ if and only if $u\in U,v\in V$
and $\{u,v\}\in E(X)$. So, it is natural to consider which bipartite
graphs have 
balanced adjacency matrices. Equivalently, which bipartite graphs have
no induced circuits of length $2\pmod 4$. We refer to such graphs as
\emph{balanced graphs}. Most results on balanced graphs are restricted
to some subclass. For instance, balanced graphs in which every induced 
circuit has length $4$ have been characterised in~\cite{Go}. It is
our aim in this paper to 
provide some insight into the structure of some additional classes of balanced
graphs. 

All graphs in this paper are connected, so, by abuse of terminology we
use the term ``graph'' to mean ``connected graph''.

We will present three main results in this paper, together with a
number of conjectures. In two of our results, we restrict our
attention to a significant family of balanced graphs, and characterise
all of the graphs in that family. First, we characterise balanced
Cayley graphs on abelian groups. Then we prove that there are no cubic
balanced planar graphs. In the remaining result, we provide a
condition on the number of vertices of a $k$-regular balanced graph.

Before we can state our characterisation of balanced Cayley graphs on
abelian groups, a number of definitions will be required.

Let $G$ be a group and $S$ a subset of $G$ such that $S$ is closed
under taking inverses, i.e. $S=S^{-1}$. The \emph{Cayley graph} of
$G$ with connection set $S$, denoted $\Cay(G,S)$, is the graph
with vertex set $G$ and edge set $\{\{g,h\}\mid gh^{-1}\in S\}$.
If $X$ and $Y$ are graphs, the lexicographic product of $X$ with
$Y$ is the graph with vertices $V(X) \times V(Y)$, where $(x,y)$
and $(x',y')$ are adjacent if and only if either $\{x,x'\} \in
E(X)$, or $x=x'$ and $\{y,y'\} \in E(Y)$. We denote by
$\overline{K_t}$  the complement of the complete graph $K_t$, i.e.
$\overline{K_t}$ has $t$ vertices and no edges. Further, we denote
by $C_l$ the cycle of length $l$ and $C_2=K_2$ for the degenerate
case. The following terminology will be 
used in our characterisation.

\begin{definition} Let $l, t $ be in $\mathbb{Z}^+$, with $l=2$ or
$l \equiv 0 \pmod{4}$ and $l\geq 8$. The lexicographic product of $C_l$ with
$\overline{K_t}$ is called an \emph{$(l,t)$-cycle}.
\end{definition}

\begin{figure}[!h]
\begin{center}
\begin{picture}(0,0)%
\includegraphics{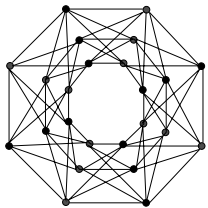}%
\end{picture}%
\setlength{\unitlength}{1105sp}%
\begingroup\makeatletter\ifx\SetFigFont\undefined%
\gdef\SetFigFont#1#2#3#4#5{%
  \reset@font\fontsize{#1}{#2pt}%
  \fontfamily{#3}\fontseries{#4}\fontshape{#5}%
  \selectfont}%
\fi\endgroup%
\begin{picture}(4486,4516)(818,-4044)
\end{picture}
\caption{$(8,3)$-cycle}\label{F:cycle}
\end{center}
\end{figure}

Figure~\ref{F:cycle} shows the example of an $(8,3)$-cycle.

We recall that the circuits of a Cayley graph $X$ are well studied
and the problem of
``understanding'' the lengths of the circuits of $X$ has attracted
considerable attention, 
see~\cite{Als} and~\cite{Wit}. For example, it is known that if
$X$ is a bipartite Cayley graph on an abelian group, then $X$ has
circuits of any even length, i.e. $4,6,\ldots,|X|$,
see~\cite{bog:pan}. Unfortunately, not much is known about the
\emph{induced} circuits of $X$. So, in this paper, we focus only
on induced circuits. Since it seems hard to get general results on
induced circuits of Cayley graphs, we restrict our interest to
balanced Cayley graphs. Specifically, in Section~\ref{S:circulant} and
Section~\ref{S:abelian}, we study balanced Cayley graphs  on
abelian groups. We will show that this class of graphs is very
restricted. Indeed, we get the following result.

\begin{theorem}\label{T:mainabelian}
If $G$ is an abelian group and $S \subseteq G$, then the graph
$\Cay(G,S)$ is balanced if and only if it is isomorphic to an
$(l,t)$-cycle.
\end{theorem}

In Section~\ref{S:divisibility}, we use a theorem of
Berge~\cite{ber:bal} and Fulkerson, Hoffman, and
Oppenheim~\cite{ful:bal} to prove the following result.

\begin{theorem}\label{T:divisibility}
If $X$ is a $k$-regular balanced graph, then the number of
vertices of $X$ is divisible by $2k$.
\end{theorem}

In Section~\ref{sec:planar}, we narrow our focus to cubic 
balanced planar graphs and we prove the following theorem.

\begin{theorem}\label{noplanar}
There is no cubic balanced planar graph.
\end{theorem}

Finally, two conjectures for balanced regular graphs are presented
in Section~\ref{S:conjectures}. These conjectures are based on the
results for balanced Cayley graphs on abelian groups, on
Theorem~\ref{noplanar} and on some exhaustive computer
computations.


\section{Balanced circulant graphs}\label{S:circulant}

It seems to the authors that the statements and the proofs in
Section~\ref{S:circulant} and~\ref{S:abelian} are neater using a
multiplicative notation. So, all groups in the following two
sections will be written using multiplicative notation. As usual, if
$G$ is a group and $a\in G$, the
symbol $\langle a\rangle$ denotes the subgroup  of $G$ generated by
$a$, and $1_G$ denotes the identity element in $G$.

Let $A_n$ denote a cyclic group of order $n$, generated by $a$,
and let $S \subseteq A_n \backslash \{1_{A_n}\}$ be closed under taking
inverses. The graph $\Cay(A_n, S)$ is said to be a \emph{circulant
graph} of order $n$. Recall that $\Cay(A_n, S)$ is bipartite if
and only if  $n$ is even and $S\subseteq\{a^i\mid i\textrm{
odd}\}$.

\begin{lemma}\label{L:complete}If $X=\Cay(A_n, S)$ is
balanced and $a,a^3 \in S$, then $X$ is $K_{n/2,n/2}$.
\end{lemma}

\noindent\textsc{Proof. } We prove, by induction on $l$, that
$a^{2l+1}$ lies in $S$ for every $l$. If $l=0$ or $1$, then there is
nothing to prove, by 
hypothesis $a,a^3$ lie in $S$. Now, assume that the claim is true for every index $i$ such that
$i\leq l-1$, and so, $\{a,a^3,\ldots,a^{2l-1}\}
\subseteq S$. Since $S=S^{-1}$, if $S=\{a,a^3,\ldots,a^{2l-1}\}$,
then $S=\{a^{i}\mid i\textrm{ odd}\}$ and there is nothing to
prove. Assume $\{a,a^3,\ldots,a^{2l-1}\}\subsetneqq S$. Let $m=2l-1$ and let $t >0$ be minimum with
$a^{m+t} \in S$. Set $t=qm +r$ where $q$ and $r$ are non-negative
integers with $0\leq r<m$. We show that $q=0$ and $r=2$ and thus
$a^{2l+1} \in S$. Note that since $S^{-1}=S$, we may assume that
$m<n/2$.

Suppose instead that $q \geq 1$. The following table lists the
vertices of an induced circuit $C$ in $X$ with length $2
\pmod{4}$. In reading the table, the following remarks may be
useful. Since $t$ is even and $m$ is odd, $q$ and $r$ have the same
parity. Each cycle ends with the vertex $a^{(q+1)m+r}=a^{m+t}$. Now,
we leave to the reader the straightforward work of 
checking that every given circuit is induced, since all differences
between non-consecutive vertices are either even, or strictly between
$m$ and $m+t$. Since $X$ is
balanced, it follows that $q=0$.

\begin{table}[!ht]
\begin{center}
\begin{tabular}{c|c|l}
$q \pmod{4}$ & $r$       & $\qquad \qquad C$ \\ \hline
$0$            & $0$         &$(1,a^m,a^{2m},\ldots, a^{(q+1)m})$\\
$0$            & $2$         &$(1,a,a^2,a^{m+2},a^{m+3},a^{m+4},a^{2m+4},a^{3m+4},\ldots,a^{qm+4}, a^{(q+1)m+r})$ \\
0              &$\geq 4$     &$(1,a,a^2,a^{m+2},a^{2m+2},\ldots,a^{(q+1)m+2},a^{(q+1)m+r-1},a^{(q+1)m+r})$ \\
$1$            & $\geq 1$         &$(1,a,a^2,a^{m+2},a^{2m+2},a^{3m+2},\ldots,a^{qm+2},a^{(q+1)m},a^{(q+1)m+r})$ \\
$2$            & 0         & $(1,a,a^2,a^{m+2},a^{2m+2},\ldots,a^{qm+2},a^{(q+1)m})$ \\
$2$            & $\geq 2$  & $(1,a^m,a^{2m},\ldots,a^{(q+1)m},a^{(q+1)m+r-1},a^{(q+1)m+r})$ \\
$3$            & $\geq 1$       & $(1,a^m,a^{2m},\ldots,a^{(q+1)m},a^{(q+1)m+r})$ \\
\end{tabular}
\end{center}
\end{table}

Since $q=0$, we have $r$ is even. If $r \geq 4$, then
$(1,a,a^2,a^{m+2},a^{m+3},a^{m+r})$ is an 
induced circuit of length $6$. Thus $r=2$ and the induction is complete. 

We have proved that
$S=A_n\backslash \langle a^2\rangle$ and so $X$ is a complete
bipartite graph.~$_\Box$

Next we consider the case that $a^3 \not \in S$.

\begin{lemma}\label{L:induction} Let $X=\Cay(A_n,S)$ be
balanced with $a \in S$, $|S|>1$ and $a^3 \not \in S$. If $l > 2$
is minimum such that $a^{l-1} \in S$, then  $l\equiv 0\pmod 4$, $l\geq
8$, $l$ divides $n$ and
$S=\{a^{il \pm 1 }\mid 0\leq i \leq (n/l)-1\}$.
\end{lemma}

\noindent\textsc{Proof. }The circuit $(1,a,a^2,\ldots,a^{l-1})$ is
induced in  $X$, therefore $l\equiv 0\pmod 4$. Moreover, since
$a^3\notin S$, we have $l\geq 8$. We use induction on $k$ to prove the
following claim, from which 
Lemma~\ref{L:induction} would follow.

\begin{sublemma}\label{C:induct}
If $k \in \mathbb{Z}^+$, then $S \cap \{a^{j}\mid -1\leq j \leq
kl+1\} = \{a^{il\pm 1} \mid 0 \leq i \leq k\}$.
\end{sublemma}

Let $k=1$. In this case, by hypothesis on $a$ and on $l$, we have
only to show that $a^{l+1}$ lies in $S$. The circuit $C =
(1,a,a^2,a^{l+1}, a^{l+2},a^{l+3},a^4,a^5,a^6,\ldots,a^{l-1})$ has
length $l+2$ and thus is not induced. Using the hypothesis that $l>2$
is minimum with $a^{l-1}\in S$, it is easy to check that
$\{1,a^{l+1}\}$, $\{1,a^{l+3}\}$, $\{a,a^{l+2}\}$, $
\{a^2,a^{l+3}\}$ are the only possible chords for $C$; hence
$a^{l+1} \in S$ or $a^{l+3} \in S$. If $a^{l+3}\in S$ and
$a^{l+1}\notin S$, then 
$(1,a^{l-1},a^l,a^{l+1},a^{l+2},a^{l+3})$ is an induced circuit of
length $6$ in $X$, a contradiction. Thus $a^{l+1}\in S$ and the
claim holds when $k=1$.

Assume that $ S \cap \{a^j\mid -1\leq j\leq kl+1\} = \{a^{il\pm 1}
\mid 0 \leq i \leq k\}$ for some $k \geq 1$. Let $t>0$ be
minimal such that $a^{kl+t+1} \in S$. We have to prove that $t=l-2$.  
If $t\equiv 0 \pmod 4$ and
$t\leq l-3$, then $(1,a^{kl+1},a^{kl+2},\ldots,a^{kl+t+1})$ is an
induced circuit of length $t+2$ in $X$, a contradiction.
Similarly, if $t\equiv 2 \pmod 4$ and $t\leq l-3$, then
$(1,a^{kl+t+1},a^{t+2},a^{t+3},\ldots,a^{l-1})$ is an induced
circuit of length $l-t$, a contradiction. This yields $t\geq l-2$.

Consider the circuit $C=(1, a^{kl -1}, a^{(k+1)l}, a^{(k+1)l+1},
a^{l+2}, a)$ of $X$. Clearly, $C$ must have a chord, so,
$a^{(k-1)l-3}\in S$, or $a^{(k+1)l-1}\in S$ 
or $a^{(k+1)l+1}\in S$. Since $(k-1)l-3\leq kl+1$ and
$(k-1)l-3\not\equiv \pm1\pmod l$, by induction
hypothesis, we have $a^{(k-1)l-3}\notin S$. This says that 
 either   $a^{(k+1)l-1}\in S$ or
$a^{(k+1)l+1}\in S$. If $a^{(k+1)l+1}\in S$ and $a^{(k+1)l-1}\notin
 S$, then the circuit 
$(1,a,a^2,\ldots,a^{l-2},a^{(k+1)l-1},a^{(k+1)l},a^{(k+1)l+1)})$
is induced of length $l+2$, a contradiction. Thus $a^{(k+1)l-1}
\in S$.

Now, to conclude the inductive argument, it remains to prove that
$a^{(k+1)l+1}$ is in $S$. Let $t>0$
be minimal such that $a^{(k+1)l-1+t} \in S$. If $t>4$, then
\[ (1,a,a^2,a^{(k+1)l+1},a^{(k+1)l+2},a^{(k+1)l+3},a^{kl+4},
a^{kl+5},a^{kl+6},\ldots,a^{(k+1)l-1})\] 
is an induced circuit of length $l+2$. Similarly, if $t=4$, then
the circuit
$(1,a^{(k+1)l-1},a^{(k+1)l},a^{(k+1)l+1},a^{(k+1)l+2},a^{(k+1)l+3})$
is induced of length $6$. This yields $t=2$, and thus $a^{(k+1)l
+1} \in S$.
 Claim~\ref{C:induct} follows.~$_\blacksquare$

Now, we are ready to conclude the proof of
Lemma~\ref{L:induction}. Since $S$ is closed under taking
inverses, Claim~\ref{C:induct} and $l\neq 4$ imply that $n-(l-1) =
kl + 1$ for some integer $k$. Thus $l$ divides $n$ and the lemma
follows.~$_\Box$

If $X$ is a graph and $u \in V(X)$, then the \emph{neighbours} of
$u$, denoted $N(u)$, are the elements of $N(u)= \{ v \in V(X) \mid \{u,v\}
\in E(X)\}$. Vertices $u,v$ of $X$ are \emph{twins} if
$N(u)=N(v)$. Being twins is an equivalence relation on $V(X)$. We say
that two vertices $u,v$ of $X$ are 
non-trivial twins if $u\neq v $ and $u,v$ are twins. If $u\in V(X)$ then we denote by $X\backslash u$ the induced subgraph of $X$ on $V(X)\setminus\{u\}$.

\begin{lemma}\label{L:twins}
If $X$ is a bipartite graph with non-trivial twins $u$ and $v$,
then $X$ is balanced if and only if $X \backslash u$ is balanced.
\end{lemma}

\noindent\textsc{Proof. } Induced subgraphs of balanced graphs are
balanced, and thus $X \backslash u$ is balanced whenever $X$ is
balanced. Conversely, suppose $X \backslash u$ is balanced. In
particular, as $N(u)=N(v)$, we have that $X\backslash v$ is also
balanced. Let $C$ be an induced circuit in $X$. If $u \not \in V(C)$,
or if  $v 
\not \in V(C)$, then $C$ is isomorphic to an induced circuit in $X
\backslash u$, or $X\backslash v$ (respectively), and thus $|V(C)| \equiv 0 \pmod 4$. If $u, v \in
V(C)$ then, since $C$ is induced and $N(u)=N(v)$, we have
$|V(C)|=4$. In either case every induced circuit of $X$ has length
$0\pmod 4$, therefore $X$ is balanced. $_\Box$

Note that the complete bipartite graph $K_{t,t}$ is a
$(2,t)$-cycle. In the following lemma, we prove that $(l,t)$-cycles
are circulant graphs.

\begin{lemma}\label{L:aux}Let $l, t$ be   in $\mathbb{Z}^+$ with $l=2$ or $l
\equiv 0 \pmod{4}$ and $l\geq 8$, $A_{lt}$ be a cyclic group of order $lt$ with
generator $a$ and $S =\{a^{il\pm 1}\mid 0\leq i\leq t-1\}$. The
circulant graph $\Cay(A_{lt}, S)$ is isomorphic to an
$(l,t)$-cycle.
\end{lemma}

\noindent\textsc{Proof. }We leave it to the reader to check that, for each $j \in \{0, \ldots,
l-1\}$, the set $\{a^{il+j}\mid 0\leq i\leq t-1\}=\langle a^l\rangle
a^j$ is the twin 
class  of $a^j$ in $\Cay(A_{lt},S)$. Furthermore, the quotient graph
through the twin equivalence relation is a cycle of length $l$. Now, the
result is straightforward. $_\Box$

As a corollary of the results we have proved, we get the
following theorem.

\begin{theorem}\label{circulant+1} A  circulant bipartite graph
$X$ with a generator in its connection set is balanced if and only
if $X$ is an $(l,t)$-cycle.
\end{theorem}

\noindent\textsc{Proof. }By
Lemmas~\ref{L:complete},~\ref{L:induction} and ~\ref{L:aux}, we
have that every balanced circulant graph $X$ with a generator in
its connection set is an $(l,t)$-cycle.  Conversely, with
Lemma~\ref{L:twins}, it is easily verified that $(l,t)$-cycles are
balanced. $_\Box$


\section{Balanced Cayley graphs on abelian groups}\label{S:abelian}

In this section  Theorem~\ref{circulant+1} is used to prove
Theorem~\ref{T:mainabelian}.

Let $G$ be an abelian group with $S \subseteq G$ such that
$X=\Cay(G,S)$ is balanced. Choose $a\in S$. Recall that if $U \subseteq
V(X)$, then $X[U]$ denotes the subgraph of $X$ induced by $U$. In
particular, we have
 $X[\langle a\rangle]=\Cay(\langle a\rangle,S\cap\langle
 a\rangle)$. Therefore $X[\langle a\rangle]$ is a circulant graph
 with a generator in its connection set.
By Theorem~\ref{circulant+1}, $X[\langle a\rangle]$ is an
$(l,t)$-cycle. Moreover, by Lemmas~\ref{L:complete}
and~\ref{L:induction}, we have that either $l=2$ or $l\equiv 0\pmod
4$, $l\geq 8$, $l$ divides $|a|$ and  
\begin{equation}\label{Eq:1}
S \cap \langle a \rangle = \{ a^{il\pm 1} \mid 0\leq i\leq
(|a|/l)-1\}.
\end{equation}
If $G=\langle a\rangle$, then
Theorem~\ref{T:mainabelian} follows. Thus we assume $G \not
=\langle a\rangle$.

For the remainder of Section~\ref{S:abelian}, $G,S,X,a$ and $l$
are fixed, with the meaning defined in the preceding paragraph. Before
going into the proof of Theorem~\ref{T:mainabelian}, we would like to
point out explicitly where the group $G$ being abelian is
used. We recall that $\Aut(X)$ contains the right regular
representation of $G$, so, if $e\in E(X)$ and $g\in G$, then
$eg=ge$ is an edge of $X$. Consequently, whenever
$H,Hb$ are two cosets of $G$ and $b\in 
S\backslash H$, then there are ``plenty'' of edges between $H$ and
$Hb$, indeed $\{h,hb\}$ is an edge for any $h\in H$. In other words,
 any edge $e$ from a vertex of $H$ to a vertex of $Hb$ determines a
 matching $\{he\mid h\in H\}$. Without the hypothesis of $G$ being
 abelian we could only say that $\{eh\mid h\in H\}$ is a family of
 edges of $X$ and no further ``structure'' on this family of edges or
 where these edges lie could be assumed.

\begin{lemma}\label{L:twocosets}
If $b \in S \backslash \langle a \rangle$, then the subgraph of $X$
induced by $\langle a \rangle \cup b \langle a \rangle$ is an
$(l,2t)$-cycle.
\end{lemma}

\noindent\textsc{Proof.}
We prove three claims, from which the lemma follows.

\begin{sublemma}\label{C:fourcycle}
If $i \equiv j \pmod l$ and $ba^i \in S$, then $ba^j\in S$.
\end{sublemma}

Assume towards a contradiction that $i\equiv j\pmod l$, $ba^i\in S$
and $ba^j\notin S$. Consider the subgraph of $X$ induced by the
union of $\{a^k\mid k \equiv 0 \pmod l\}$ and $\{ba^k \mid k \equiv i
\pmod l\}$.  Note that each vertex in this subgraph has degree
equal to $| S \cap \{ba^k \mid k \equiv i \pmod l\}|$. In particular,
since $1_G$ and $a^{i-j}$ have the same degree and $\{1_G, ba^i\}$
is an edge but $\{a^{i-j}, ba^i\}$ is not, there exists $k \equiv
i \pmod l$ such that $\{a^{i-j}, ba^k\}$ is an edge but $\{1_G,
ba^k\}$ is not. Now, $j+1-k\equiv 1\pmod l$, so, by Eq.~\ref{Eq:1}, we have
$a^{j+1-k}\in S$. Therefore, $(1_G,ba^i,ba^{j+1},ba^k,a^{i-j},a) $ is an
induced  circuit  in $X$ and it has
length $6$, a contradiction. Claim~\ref{C:fourcycle}
follows.~$_\blacksquare$

Suppose $l=2$, i.e. $X[\langle a\rangle]$ is a complete bipartite graph.
Then $X[\langle a\rangle]$ is a $(2,t)$-cycle where $2t=|a|$, and
$a^i \in S$ for all odd $i$. Since $b\in S$,
Claim~\ref{C:fourcycle} yields $ba^j\in S$ for all even $j$. This
yields 
\begin{equation}\label{Eq:2}
S\cap (\langle a\rangle\cup b\langle a\rangle)=\{a^i,ba^j\mid
i\equiv 1\pmod 2, j\equiv 0\pmod 2\}.
\end{equation}

We leave it to the reader to prove that  Eq.~\ref{Eq:2} yields 
$X[\langle a \rangle \cup b \langle a \rangle]$ is a
$(2,2t)$-cycle. 

Hence we assume that $l\geq 8$, i.e. $X[\langle a\rangle]$ is not
a complete bipartite graph.

\begin{sublemma}\label{C:edge+-2}
If $ba^i \in S$, then $ba^{i+2} \in S$ or $ba^{i-2}\in S$.
\end{sublemma}

We argue by contradiction, so, assume $ba^{i-2},ba^{i+2}\notin S$. The
circuit $C=(1_G, a, a^2, ba^{i+2}, ba^{i+3}, \ldots, 
ba^{i+l-1}, ba^i)$ has length $l+2$. So $C$ is not induced. Using
Eq.~\ref{Eq:1} to study the possible chords
in $C$ and applying Claim~\ref{C:fourcycle}, we get that either
$ba^{i-2}\in S$ or  $ba^{i+m} \in S$ for some $2 \leq m \leq l-4$.
Since $ba^{i-2},ba^{i+2}\notin S$, we have $ba^{i+m}\in S$ for some
$4\leq m\leq l-4$. Let $m > 0$ 
be minimum with
$ba^{i+m} \in S$. The circuit $(1_G, ba^i, ba^{i+1}, \ldots,
ba^{i+m})$ is induced in $X$ of length $m+2$. Therefore $m \equiv 2
\pmod 4$. Hence  $6 \leq m \leq l-6$. The $10$~-~cycle
$(1_G, ba^i, ba^{i+1}, ba^{i+2}, ba^{i+3}, a^3, ba^{i+m+3}, ba^{
i+m+2}, ba^{i+m+1}, ba^{i+m})$ has only one possible chord, namely
$\{1_G,ba^{i+m+2}\}$. Since $X$ is balanced, we have $ba^{i+m+2}\in
S$. Now, the circuit $(1,ba^{i+m+2},a^2,ba^{i+2},ba^{i+1},ba^i)$ is
induced in $X$ and has length $6$, a contradiction.~$_\blacksquare$

Since $b \in S$, Claim~\ref{C:edge+-2} says that either $ba^2\in
S$ or $ba^{-2}\in S$. Since the roles of $a$ and $a^{-1}$ are
interchangeable in our arguments, we assume, throughout the rest of
Lemma~\ref{L:twocosets},  that $ba^2\in S$.

\begin{sublemma}\label{C:onlyedges}
If $ba^i \in S$, then $i \equiv 0 \pmod l$ or $i \equiv 2 \pmod
l$.
\end{sublemma}

We argue by contradiction. So, let $i\geq 0$ such that $ba^i\in S$
and $i\not\equiv 0,2\pmod l$, in particular, pick such an $i$ as
small as possible. Claim~\ref{C:fourcycle} yields $4\leq i\leq
l-2$. Write $i=m+2$, so, $2\leq m\leq l-4$. Since the circuit
$(1_G, ba^2,ba^3,\ldots,ba^{2+m})$ is induced, we have $m \equiv 2
\pmod 4$, and thus $2\leq m\leq l-6$.

The circuit $C=(1_G,ba^{2+m},ba^{3+m},\ldots,ba^{l-1},b)$ of $X$ has
length $l-m\equiv 2\pmod 4$. Therefore, $C$ has a chord. It follows that there
exists $m'$ with $ba^{2+m+m'} \in S$ and $4+m\leq 2+m+m'\leq l-2$.
Pick $m'$ minimal with these properties.

By Claim~\ref{C:edge+-2}, either $ba^{m+m'}\in S$ or
$ba^{4+m+m'}\in S$. If $ba^{m+m'} \in S$, then $m'=2$ by
minimality of $m'$. However, when
 $m'=2$, the circuit $(1_G,ba^2,ba^3,a^3,ba^{5+m}$, $ba^{4+m})$ is induced.
This contradiction implies that $ba^{4+m+m'}\in S$ and $m'>2$.

Claim~\ref{C:edge+-2} with $i=2+m$ implies that either $ba^m \in S$
or $ba^{4+m}\in S$.
Since $m' > 2$, the minimality of $m'$ implies that $ba^{4+m} \not \in S$.
Thus $ba^m \in S$. If $m>2$, this contradicts the minimality of $i$,
so $m=2$. In particular, $ba^{6+m'}\in S$.

Consider the circuit
$C=(1_G,ba^2,ba^3,a^3,ba^{7+m'},ba^{6+m'})$. Since $m=2$, we have $2
\le m' \le l-6$. By Eq.~\ref{Eq:1},
the only possible chord in $C$ is $\{ba^2,ba^{7+m'}\}$ and thus
$a^{m'+5}\in S$. Since  $m'+5\leq l-1$, Eq.~\ref{Eq:1}
implies that $m'+5=l-1$.
Since $ba^{2+m+m'}\in S$, we have $ba^{l-2}\in S$.
By Claim~\ref{C:fourcycle}, $ba^{-2}\in S$ and thus the circuit
$(1,ba^2,ba,a^3,ba^5,ba^4)$ is induced in $X$.
This contradiction finally implies Claim~\ref{C:onlyedges}.~$_\blacksquare$

Since $b,ba^2\in S$, Claims~\ref{C:fourcycle} and~\ref{C:onlyedges} show that
 $ba^i\in S$ if and only if $i\equiv 0,2\pmod l$. Thus 
\begin{equation}\label{Eq:3}
S \cap ( \langle a
\rangle \cup b \langle a \rangle ) = \{a^i, ba^j \mid i \equiv \pm 1
\pmod l, j \equiv 0,2 \pmod l\}.
\end{equation}
We leave it to the reader to check that Eq.~\ref{Eq:3} yields that $X[\langle
  a\rangle\cup b\langle a\rangle]$ is an $(l,2t)$-cycle (show that
  $a^i$ and $ba^j$  
are twins in $X[ \langle a \rangle \cup b \langle a \rangle]$ if
and only if $j \equiv i + 1\pmod l$). The proof of
  Lemma~\ref{L:twocosets} is complete. $_\Box$  

\noindent\textsc{Proof of Theorem \ref{T:mainabelian}. } Let
$\widehat{X}$ denote the quotient graph $X/\langle a\rangle$. That
is, $\widehat{X}$ has a vertex for each coset of $\langle a
\rangle$, and the vertices $b_1\langle a \rangle$ and $b_2\langle
a \rangle$ are adjacent in $\widehat{X}$ if and only if
$b_1\langle a \rangle\cap b_2\langle a \rangle =\emptyset$ and
there exists $u_1 \in b_1\langle a \rangle$ and $u_2 \in
b_2\langle a \rangle$  with $\{u_1,u_2\} \in E(X)$.

\noindent\textsc{Case 1: }$l=2$, i.e. $X[\langle a\rangle]$ is a
complete bipartite graph.

\noindent Let $(\alpha_1,\alpha_2,\ldots,\alpha_k)$ be an induced
circuit in $\widehat{X}$. Note that by Lemma~\ref{L:twocosets},
$X[\alpha_i\cup\alpha_{i+1}]$ is a complete bipartite graph
(indices $\mod k$). If $k\geq 4$ then it is straightforward to
exhibit an induced circuit of $X$ of  length $2\pmod 4$. For
instance, if $k\equiv 0\pmod 4$ then pick a path of length $2$ in
$X[\alpha_1]$ and a path of length $2$ in  $X[\alpha_{k-1}]$ and
extend it to an induced circuit of $X$ using a unique vertex from
each $X[\alpha_i],$ where $ i \neq 1,{k-1}$.

\vspace{.7cm}
\rput(1.7,0){\psellipse(0,0)(.5,1)\psellipse(2,0)(.5,1)\psellipse(4,0)(.5,1)\psellipse(7,0)(.5,1)\psellipse(9,0)(.5,1)\psellipse(11,0)(.5,1)
\rput(0,1.2){$\alpha_1$}\rput(2,1.2){$\alpha_2$}\rput(4,1.2){$\alpha_3$}\rput(7,1.2){$\alpha_{k-2}$}\rput(9,1.2){$\alpha_{k-1}$}\rput(11,1.2){$\alpha_k$}
\cnode(0,.5){.1}{1}\cnode*(0,-.5){.1}{2}\cnode(2,0){.1}{3}\cnode*(4,0){.1}{4}\cnode(7,0){.1}{5}\cnode*(9,.5){.1}{6}\cnode(9,-.5){.1}{7}\cnode*(11,.5){.1}{8}
\cnode(5,0){0}{a}\cnode(6,0){0}{b}
\ncline{1}{2}\ncline{2}{3}\ncline{3}{4}\ncline[linestyle=dotted]{4}{a}\ncline[linestyle=dotted]{b}{5}\ncline{5}{6}\ncline{5}{6}\ncline{6}{7}\ncline{7}{8}
\ncangle[angleA=90,armA=.2cm, armB=.2cm,angleB=90]{1}{8}}
\vspace{.7cm}

All other cases are fairly similar. This proves that every induced
circuit in $\widehat{X}$ has length  $3$, and thus $\widehat{X}$
is complete. Now, it is easy to check, using Lemma~\ref{L:twocosets}
and the fact that
$\widehat{X}$ is complete, that $X$ is
isomorphic to a $(2,t|G|/|a|)$-cycle.

\noindent\textsc{Case 2: }$l\geq 8$.

We first prove a claim.
\begin{sublemma}\label{sbL:1}The graph $\widehat{X}$ is
  complete. Further, let $\alpha_1,\alpha_2,\alpha_3$ be  three distinct
  vertices of $\widehat{X}$. If $x_2\in V(X[\alpha_2])$  has twins
  $x_1\in V(X[\alpha_1])$ in $X[\alpha_1\cup\alpha_2]$ and $x_3\in
  V(X[\alpha_3])$ in $X[\alpha_2\cup \alpha_3]$, then $x_1$ and
  $x_3$ are twins in $X[\alpha_1\cup \alpha_3]$. 
\end{sublemma} 

Let $(\alpha_1,\alpha_2,\alpha_3)$
be a path of length $3$ in $\widehat{X}$. For  $i\in\{1,2,3\}$, choose $x_i \in V(X[\alpha_i])$ such that
$x_1$ and $x_2$ are twins in $X[\alpha_1\cup\alpha_2]$ and $x_2$
and $x_3$ are twins in  $X[\alpha_2\cup\alpha_3]$ (this is feasible
because, by Lemma~\ref{L:twocosets}, the graphs $X[\alpha_1\cup \alpha_2]$
and $X[\alpha_2\cup \alpha_3]$ are $(l,2t)$-cycles). Since the
multiplication on the right by an element of $\langle a\rangle$ is
an automorphism of $X[\alpha_1\cup \alpha_2]$, we have that 
 $x_1a^j$ and $x_2a^j$ are twins in $X[\alpha_1\cup\alpha_2]$
for any $j$. Similarly, $x_2a^j$ and $x_3a^j$ are twins in
$X[\alpha_2\cup\alpha_3]$. This yields that
$(x_1,x_1a,x_2a^2,x_3a,x_3,x_2a^{-1})$ is a circuit of length $6$
in $X$, therefore it has a chord. Since $l\geq 8$, we have that
$\{x_2a^2,x_2a^{-1}\} \not \in E(X)$, therefore either
$\{x_1,x_3a\}$ or $\{x_1a,x_3\}$ is an edge of $X$. Without loss of
generality, we may assume that $\{x_1,x_3a\}$ is an edge (the role of
$a$ and $a^{-1}$ is interchangeable). In particular,
$\{\alpha_1, \alpha_3\} \in E(\widehat{X})$, and thus
$\widehat{X}$ is complete.

Now, it remains to prove that $x_1$ and $x_3$ are twins in
$X[\alpha_1\cup\alpha_3]$. Since $\widehat{X}$ is complete,
Lemma~\ref{L:twocosets} yields that $X[\alpha_1\cup\alpha_3]$ is an
$(l,2t)$-cycle. Moreover, since $\{x_1,x_3a\} \in E(X)$, we have that 
either $x_1$ is a twin with $x_3$ in
$X[\alpha_1\cup\alpha_3]$ or $x_1$ is a twin with $x_3a^2$ in
$X[\alpha_1\cup\alpha_3]$. In the latter case,
$(x_1,x_1a,x_1a^2,\ldots,x_1a^{l-5},x_3a^{l-2},x_2a^{l-1})$ is an
induced circuit 
of length $l-2$, and by this contradiction, $x_1$ and $x_3$
are twins in $X[\alpha_1\cup\alpha_3]$. The claim is proved.~$_\blacksquare$

Now, we leave it to the reader to check that Lemma~\ref{L:twocosets}
and Claim~\ref{sbL:1} yield that $X$ is a $(l,t|G|/|a|)$-cycle. The
proof of Theorem~\ref{T:mainabelian} is complete.~$_\Box$

\section{The number of vertices in a regular balanced
graph}\label{S:divisibility}

The graph $X=\Cay(G,S)$ is regular with degree $|S|$.
Theorem~\ref{T:mainabelian} implies that if $X$ is balanced and
$G$ is abelian, then $|G|$ is divisible by $2|S|$. We prove
Theorem~\ref{T:divisibility}, which states  that this divisibility
criterion holds for all regular balanced graphs. To do this we
need first some terminology and a well-known result in linear
programming on balanced matrices.

We recall that if $A$ is an $n\times m$ $\{0,1\}$-matrix, then the
set partitioning polytope defined by $A$ is the set $R(A)=\{x\in
\mathbb{R}^m\mid Ax=\textbf{1}_n,\textbf{0}_m\leq x\leq
\textbf{1}_m\}$, where $\mathbf{1}_n$, respectively
$\mathbf{1}_m$, denotes a column vector of length $n$,
respectively $m$, whose entries are all equal to $1$ and
$\mathbf{0}_m$ denotes the zero vector of length $m$. Note that $R(A)$
is a convex polytope. A set
partitioning polytope is said to be \emph{integral} if all its
vertices (i.e. extremal points) have only integer-valued
components.

The following characterization of balanced matrices is due to
Berge~\cite{ber:bal} and Fulkerson, Hoffman, and
Oppenheim~\cite{ful:bal}.

\begin{theorem}[Berge~\cite{ber:bal}, Fulkerson, Hoffman, and
Oppenheim~\cite{ful:bal}]\label{T:bergefulkerson}
If $A$ is  a balanced matrix, then 
$R(A)$ is integral.
\end{theorem}

We note that a proof of Theorem~\ref{T:bergefulkerson} is also
in~\cite{con:fag} (see Theorem~$2.1$), where a more general result
is proved. Now, Theorem~\ref{T:divisibility} is a corollary of
Theorem~\ref{T:bergefulkerson}.

\noindent\textsc{Proof of Theorem~\ref{T:divisibility}. }Let $X$
be a $k$-regular balanced graph and $A$ be the bipartite adjacency
matrix of $X$. Since $X$ is regular, $A$ is a square matrix of
size $n=\frac{1}{2}|V(X)|$. Since the vector
$\frac{1}{k}\textbf{1}_n$ lies in the set
partitioning polytope $R(A)$, we have $R(A)\neq \emptyset$. Thus,
by Theorem~\ref{T:bergefulkerson}, there exists a vertex $x$ of $R(A)$
with integer-valued components. As $x$ lies in $R(A)$, we have
$\textbf{0}_n\leq x\leq\textbf{1}_n$, so, $x$  has components in
$\{0,1\}$. Set $t=\sum_{i=1}^n x_i$. The 
equation $Ax=\mathbf{1}_n$ and the fact that $X$ is $k$-regular
yield $tk =n=1/2|V(X)|$. The theorem follows.~$_\Box$

\section{Cubic balanced planar graphs}\label{sec:planar}

In this section we deal with cubic planar graphs and we prove
Theorem~\ref{noplanar}.

\begin{center}
\begin{figure}[!h]
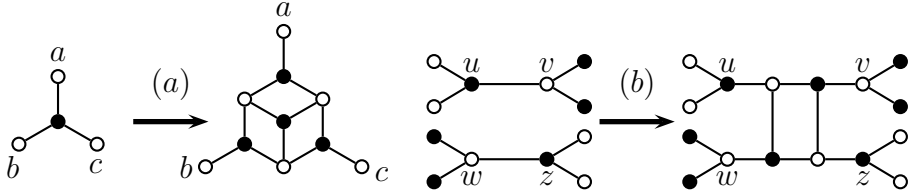

\rput(1.7,.8){ \cnode*(0,0){.1}{1}\cnode(0,.6){.1}{2}
\cnode(.5196,-.3){.1}{3}\cnode(-.5196,-.3){.1}{4}
\ncline{1}{2}\ncline{1}{3}\ncline{1}{4}
\cnode*(3,0){.1}{5}\cnode(3,-.6){.1}{6}
\cnode(3.5196,.3){.1}{7}\cnode(2.4804,.3){.1}{8}
\cnode*(3.5196,-.3){.1}{9}\cnode*(2.4804,-.3){.1}{10}
\cnode*(3,.6){.1}{11}\cnode(3,1.2){.1}{12}\cnode(1.9608,-.6){.1}{13}
\cnode(4.0392,-.6){.1}{14}
\ncline{5}{6}\ncline{5}{7}\ncline{5}{8}\ncline{7}{9}
\ncline{8}{10}\ncline{9}{6}\ncline{10}{6}\ncline{7}{11}\ncline{8}{11}
\ncline{11}{12}\ncline{13}{10}\ncline{14}{9}
\psline[linewidth=2pt]{->}(1,0)(2,0)\rput(1.5,.5){$(a)$}}
\rput(7.2,1.3){
\cnode*(0,0){.1}{21}\cnode(1,0){.1}{22}\cnode(0,-1){.1}{23}
\cnode*(1,-1){.1}{24}
\cnode(-.5,.3){.1}{25}\cnode*(1.5,.3){.1}{26}\cnode(-.5,-.3){.1}{27}
\cnode*(1.5,-.3){.1}{28}\cnode*(-.5,-.7){.1}{29}\cnode*(-.5,-1.3){.1}{30}
\cnode(1.5,-.7){.1}{31}\cnode(1.5,-1.3){.1}{32}
\ncline{21}{22}\ncline{23}{24}\ncline{21}{25}\ncline{21}{27}\ncline{22}{26}
\ncline{22}{28}\ncline{29}{23}\ncline{30}{23}\ncline{31}{24}\ncline{32}{24}
\psline[linewidth=2pt]{->}(1.7,-.5)(2.7,-.5)\rput(2.2,0){$(b)$}
\rput(0,.25){$u$}\rput(1,.25){$v$}\rput(0,-1.25){$w$}\rput(1,-1.25){$z$}
\cnode(2.9,.3){.1}{41}\cnode(2.9,-.3){.1}{42}\cnode*(3.4,0){.1}{43}
\cnode*(2.9,-.7){.1}{44}\cnode*(2.9,-1.3){.1}{45}\cnode(3.4,-1){.1}{46}
\cnode(4,0){.1}{47}\cnode*(4,-1){.1}{48}\cnode*(4.6,0){.1}{49}
\cnode(4.6,-1){.1}{50}\cnode(5.2,0){.1}{51}\cnode*(5.2,-1){.1}{52}
\cnode*(5.7,.3){.1}{53}\cnode*(5.7,-.3){.1}{54}
\cnode(5.7,-.7){.1}{55}\cnode(5.7,-1.3){.1}{56}
\ncline{41}{43}\ncline{42}{43}\ncline{46}{44}\ncline{45}{46}\ncline{46}{48}
\ncline{47}{43}\ncline{47}{48}\ncline{48}{50}\ncline{47}{49}\ncline{49}{50}
\ncline{49}{51}\ncline{50}{52}\ncline{51}{53}\ncline{51}{54}\ncline{55}{52}
\ncline{56}{52}
\rput(-2.5,1){$a$}\rput(-5.5,.4){$a$}\rput(-6.1,-1.1){$b$}\rput(-5,-1.1){$c$}
\rput(-3.8,-1.1){$b$}\rput(-1.2,-1.2){$c$}
\rput(3.4,.25){$u$}\rput(3.4,-1.25){$w$}\rput(5.2,.25){$v$}
\rput(5.2,-1.25){$z$}} 
\caption{$(a)$ diamond inflation, $(b)$
$A_1$ subdivision}\label{fig1}
\end{figure}
\end{center}

Batagelj~\cite{Ba} proved that all $3$-connected cubic bipartite
planar graphs can be obtained from the cube by a succession of two
elementary operations. The first operation is called \emph{diamond
inflation} of a vertex and it replaces a vertex with a
``diamond'', see Figure~\ref{fig1}$(a)$. The second operation is
called $A_1$ \emph{subdivision} and it applies to a pair of
non-adjacent edges $\{u,v\}$, $\{w,z\}$, see
Figure~\ref{fig1}$(b)$.
\vspace{2cm}
\vspace{-1cm}

We state Batagelj's theorem precisely.
\begin{theorem}[Batagelj~\cite{Ba}]\label{T:Batagelj}
Every $3$-connected cubic bipartite planar graph can be obtained
from the cube  by a succession of diamond inflations and
$A_1$ subdivisions. The operations can be chosen such that the
intermediate graphs are planar and bipartite.
\end{theorem}

\begin{lemma}\label{Lemma:1}Let $X$  be a cubic bipartite planar graph with
  connectivity 
$\kappa(X)=2$. There exists a $2$-vertex cut $\{u,v\}$ such that
$X\backslash\{u,v\}$ has a component $Y$ with the following
property: there exist two nonadjacent vertices $a,b$ in $Y$ such that
$\overline{Y}=(V(Y),E(Y)\cup \{ab\})$ is a cubic $3$-connected bipartite
planar graph.   
\end{lemma}

\textsc{Proof. }Choose a 2-vertex cut-set $\{a,b\}$ so as to minimise
 the order of one of the resulting components, $C$.  If $a$ has a
 single neighbour $a'$ in $C$, then $\{a', b\}$ is a
 2-vertex cut-set that contradicts the  minimality of $C$; similarly
 for $b$.  So each of $a$ and $b$ has two  neighbours in $C$.  A
 counting argument, using the regularity and the 
 bipartition of the graph $X$, forces $a$ and $b$ to have opposite
 colours, so the neighbours of $a$ and $b$ are all distinct.  Further,
 $a$ and $b$ are not 
 adjacent, since there must be another component to which at least one of them
 is joined by an edge, and each has only three neighbours.

 We claim that the choice of $Y$ to be the induced subgraph on $V(C) \cup
 \{a,b\}$ (with the neighbours of $a$ and $b$ that are not in $C$ as the
 2-vertex cut-set) satisfies the claims of this lemma.  The minimality of $C$
 is sufficient to ensure that $\overline{Y}$ is 3-connected, and it is clearly
 cubic, bipartite and planar.~$_\Box$

\begin{center}
\begin{figure}[!h]
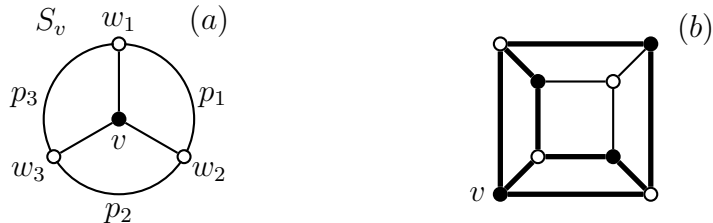

\rput(4,0){
\cnode*(0,0){.1}{v}
\cnode(0,1){.1}{w_1}
\cnode(.866,-.5){.1}{w_2}
\cnode(-.866,-.5){.1}{w_3}
\ncline{v}{w_1}\ncline{v}{w_2}\ncline{v}{w_3}
\psarc(0,0){1}{94}{206}
\psarc(0,0){1}{214}{326}
\psarc(0,0){1}{334}{86}
\rput(-.9,1.3){$S_v$}\rput(1.25,.3){$p_1$}
\rput(-1.25,.3){$p_3$}\rput(0,-1.25){$p_2$}
\rput(0,1.3){$w_1$}\rput(1.2,-.7){$w_2$}\rput(-1.2,-.7){$w_3$}
\rput(1.2,1.3){$(a)$}\rput(0,-.3){$v$}
\rput(6,0){
\rput(1.6,1.2){$(b)$}
\cnode*(1,1){.1}{1}
\cnode(-1,1){.1}{2}
\cnode*(-1,-1){.1}{3}
\cnode(1,-1){.1}{4}
\cnode(.5,.5){.1}{5}
\cnode*(-.5,.5){.1}{6}
\cnode(-.5,-.5){.1}{7}
\cnode*(.5,-.5){.1}{8}
\ncline[linewidth=2pt]{1}{2}\ncline[linewidth=2pt]{2}{3}
\ncline[linewidth=2pt]{3}{4}\ncline[linewidth=2pt]{4}{1}
\ncline{5}{6}\ncline[linewidth=2pt]{6}{7}
\ncline[linewidth=2pt]{7}{8}\ncline{8}{5}
\ncline{1}{5}\ncline[linewidth=2pt]{2}{6}
\ncline[linewidth=2pt]{3}{7}\ncline[linewidth=2pt]{4}{8}
\rput(-1.3,-1){$v$}}}
\vspace{.9cm}
\caption{The subgraph $S_v$}\label{Fig:2}
\end{figure}
\end{center}

Let $X$ be a $3$-connected cubic bipartite planar
 graph, $v$ be a vertex of $X$ and $w_1,w_2,w_3$ be the
  three neighbours of $v$. Let
  $p_i=(w_i,w_{i1},\ldots,w_{i{n_i}},w_{i+1})$ be the path in $X$ from
  $w_i$ to $w_{i+1}$ (indices  $\mod 3$) such that the circuit $(v,p_i)$
  is the boundary of a face of $X$, see Fig.~\ref{Fig:2}. Let us
  denote by $S_v$ the subgraph of $X$ with vertex set
  $\{v,w_1,w_2,w_3,w_{ij}\mid 1\leq i\leq 3,1\leq j\leq n_i\}$ and
  edge set
  $\{vw_1,vw_2,vw_3,w_iw_{i1},w_{in_i}w_{i+1},w_{ij}w_{i(j+1)}\mid
  1\leq i\leq 3,1\leq j\leq n_i-1\}$, i.e. $S_v$ is the graph in
  Fig.~\ref{Fig:2}$(a)$. For instance, the bold edges in
  Fig.~\ref{Fig:2}$(b)$ show the subgraph $S_v$ for the cube.

\vspace{1.5cm}

\begin{sublemma}\label{Sublemma:1}
$S_v$ is an induced subgraph of $X$.
\end{sublemma}

\noindent\textsc{Proof. }This is true for the cube. Now use
Theorem~\ref{T:Batagelj} and induction.~$_\blacksquare$

\begin{sublemma}\label{Sublemma:2}
$\cap_{v\in V(X)}E(S_v)=\emptyset$.
\end{sublemma}

\noindent\textsc{Proof. }This is readily true for the cube. Now use
Theorem~\ref{T:Batagelj} and induction to conclude.~$_\blacksquare$

If $X=(V,E)$ is a graph and $e\in E$, then we denote by $X\backslash
e$ the graph $(V,E\backslash\{e\})$.
\begin{lemma}\label{Lemma:2}
Let $X$ be a $3$-connected cubic bipartite planar graph. Then $X$ is
unbalanced. If $e$ is an edge of $X$, then $X\backslash e$ is
unbalanced.
\end{lemma}

\noindent\textsc{Proof. }Let
$v$ be  a vertex of a $3$-connected cubic bipartite planar graph
$X$. We claim that $S_v$ is unbalanced. If the path
$p_i$ has length $0\pmod 4$, for some $i\in\{1,2,3\}$, then the
circuit $(v,p_i)$ is induced and 
has length $2\pmod 4$, so $S_v$ is unbalanced. Finally, if $p_i$ has
length $2\pmod 4$, for every $i$, then the circuit $(p_1,p_2,p_3)$
is induced in $S_v$ and has length $2\pmod 4$, so $S_v$ is
unbalanced. Thus our claim is proved.
 
Now, by Claim~\ref{Sublemma:1}, $S_v$  is an induced subgraph of
$X$. Therefore, $X$ is unbalanced. 

Let $e$ be an edge of $X$. By Claim~\ref{Sublemma:2}, there exists
$v\in V(X)$ such that $e\notin E(S_v)$. Now, by
Claim~\ref{Sublemma:1}, $S_v$ is an induced subgraph of $X\backslash
e$. Therefore, $X\backslash e$ is unbalanced.~$_\Box$ 

\noindent\textsc{Proof of Theorem~\ref{noplanar}. }Let $X$ be a cubic
bipartite planar graph. We have to prove that $X$ is unbalanced. By
Lemma~\ref{Lemma:2}, we may assume that $\kappa(X)=2$. 
By Lemma~\ref{Lemma:1}, there exist nonadjacent $a,b$ in $V(X)$
and $Y$ an induced subgraph of $X$ such that
$\overline{Y}=(V(Y),E(Y)\cup \{e\})$ is a $3$-connected bipartite planar graph,
where $e=ab$. Now, by Lemma~\ref{Lemma:2},
$Y=\overline{Y}\backslash e$ is
unbalanced. So, $X$ is unbalanced. The proof of
Theorem~\ref{noplanar} is 
complete.~$_\Box$


\section{Conjectures}\label{S:conjectures}

The conjectures presented here are supported by computer searches performed with the invaluable help of GAP \cite{GAP}, including an exhaustive analysis of the graphs from Gordon Royle's web page.

\begin{conjecture}\label{C:vtg}
If $X$ is a connected vertex-transitive balanced graph, then $X$ is an $(l,t)$-cycle.
\end{conjecture}

The truth of Conjecture~\ref{C:vtg} would imply that the only
connected vertex-transitive balanced graphs of odd degree are the
complete bipartite graphs.

Every cubic balanced graph known to the authors has non-trivial
twins. So, we present the following conjecture which has been
verified for all graphs with fewer than $54$ vertices.

\begin{conjecture}\label{C:twins}
If $X$ is a cubic balanced graph, then $X$ has non-trivial twins.
\end{conjecture}

The truth of Conjecture~\ref{C:twins} might shed some new light
on the graph structure of a cubic balanced graph. For example the
following conjectures, interesting in their own right, are implied
by the validity of Conjecture~\ref{C:twins}.

\begin{itemize}
\item[$(i)$] $K_{3,3}$ is the only connected vertex-transitive cubic
  balanced graph; 
\item[$(ii)$] every cubic balanced graph has girth four; 
\item[$(iii)$] Conforti-Rao conjecture is true for cubic balanced
  graphs, see \cite{cor:com} page $54$. 
\end{itemize}

We leave it as an exercise to the reader to prove that
Conjecture~\ref{C:twins} yields $(i),(ii)$ and $(iii)$. 
As a matter of curiosity we point out that the truth of
Conjecture~\ref{C:twins} would also yield Theorem~\ref{noplanar}. Indeed, it
is easy to show that cubic bipartite graphs with non-trivial twins
are not planar graphs.

We conclude by reporting the number of ``small'' cubic balanced
graphs; $f(d)$ denotes the number of connected cubic balanced
graphs on $d$ vertices.\footnote{for the graphs corresponding to
the values of $d$ contact the second author}
\begin{center}
\begin{tabular}{|c|cccccc|}\hline
$d$  & 6&12&18&24&30&36\\ \hline
$f(d)$& 1&1&4&13&74&527\\ \hline
\end{tabular}
\end{center}



\begin{thebibliography}{6}
\bibitem{Als}B.Alspach, Lifting Hamilton cycles of quotient
graphs, \emph{Discrete Mathematics} \textbf{78}, 25--36, 1989.

\bibitem{Ba}V.Batagelj, Inductive definition of two restricted
classes of triangulations, \emph{Discrete Mathematics}
\textbf{52}, 113--121, 1984.

\bibitem{ber:bal}
C.Berge, Balanced matrices, {\em Mathematical Programming}
\textbf{2}, 19--31, 1972.

\bibitem{bog:pan}
Z.R.Bogdanowicz, Pancyclicity of connected circulant graphs, {\em
J. Graph Theory} \textbf{22}(2), 167--174, 1996.

\bibitem{con:fag}M.Conforti, G.Cornu\'{e}jols and
K.Vu\v{s}kovi\'{c}, Balanced Matrices, \emph{Discrete Applied
Mathematics} \textbf{306}, 2411-2437, 2006.

\bibitem{cor:com}
G.Cornu{\'e}jols, {\em Combinatorial Optimization: Packing and
Covering}, Volume~\textbf{74} of
  {\em CBMS-NSF Regional Conference Series in Applied Mathematics},
  Society for Industrial and Applied Mathematics (SIAM), Philadelphia,
  PA, 2001.

\bibitem{ful:bal}
D.R.Fulkerson, A.J.Hoffman and  R.Oppenheim, On balanced matrices,
{\em Mathematical Programming Study} \textbf{1}, 120--132, 1974.

\bibitem{GAP}
The~GAP Group, Gap -- groups, algorithms, and programming, version
4.4, 2005, \verb+(http://www.gap-system.org)+.

\bibitem{Go}G.Goss, Perfect elimination and chordal bipartite graphs,
  \emph{Journal of Graph Theory} \textbf{2}, 155--163, 1978.

\bibitem{Ta}A.Tamir, A class of balanced matrices arising from
  location problems, \emph{SIAM Journal on Algebraic and Discrete
  Methods} \textbf{4}, 363--370, 1983.  

\bibitem{Wit}D.Witte and J.A.Gallian, A survey: Hamilton cycles in
Cayley graphs, \emph{Discrete Mathematics} \textbf{51}, 293--304,
1984.
\end{thebibliography}
\end{document}